\documentclass[12pt]{amsart}

\usepackage{hyperref}

\usepackage{amsmath,amssymb,epsfig,amscd}

\newcommand{\new}{\newcommand}

\new{\bracket}[1]{\langle #1
\rangle} 
\new{\defequals}{\stackrel{\operatorname{def}}{=}}
\new{\intprod}{\!:\!}
\new{\dotprod}{\!\cdot\!}
\new{\into}{Hookrightarrow} 
\new{\tensor}{\otimes}
\new{\tnsr}{\otimes} 
\new{\iso}{\cong} 
\new{\union}{\cup}
\new{\maps}{\colon} 
\new{\goesto}{\rightarrow}
\new{\BI}{\operatorname{BI}}
\new{\Vol}{\operatorname{Vol}}

\newtheorem{proposition}{Proposition}[section]
\newtheorem{theorem}{Theorem}[section]

\newtheorem{lemma}{Lemma}[section]
\newtheorem{corollary}{Corollary}[section]

\newcounter{letter}

\newenvironment{pf}{
\begin{trivlist} \item[Proof:]}{\qed \end{trivlist}}

\new{\lieg}{{\mathfrak{g}}}               
\new{\lieh}{{\mathfrak{h}}}               
\new{\liek}{{\mathfrak{k}}}               

\new{\RR}{{\mathbb{R}}}
\new{\NN}{{\mathbb{N}}}
\new{\CC}{{\mathbb{C}}}
\new{\ZZ}{{\mathbb{Z}}}

\new{\Lie}{\operatorname{Lie}} 
\new{\Ad}{\operatorname{Ad}}
\new{\ad}{\operatorname{ad}}
\new{\vol}{\operatorname{vol}}
\new{\Cartan}{\operatorname{Cartan}}
\new{\Kirwan}{\operatorname{Kirwan}}
\new{\Project}{\operatorname{Project}}
\new{\Stab}{\operatorname{Stab}}
\new{\dint}{{\mathbf{d}}}
\new{\compact}{\text{compact}}

\begin{document}

\title{Witten's Nonabelian Localization for Noncompact Hamiltonian Spaces}

\author{Stephen F. Sawin}
        \address{ Department of Math and C. S. \\ Fairfield University \\ Fairfield, CT
06825-5195}
    \email{ sawin@cs.fairfield.edu} 
 
\begin{abstract}
For a finite-dimensional (but possibly noncompact) symplectic manifold with a compact group
acting with a proper moment map, we show that the square of the moment
map is an equivariantly perfect Morse function in the sense of Kirwan, and that the set of critical
points of the square of the moment map is a countable discrete union of compact
sets.  We show that  certain  integrals of  equivariant cohomology
classes localizes as a sum of contributions from these compact
critical sets, and we bound the contribution from each critical set.  In the case  (1) that the contribution from higher critical sets grows slowly enough that the overall integral converges rapidly and (2) that $0$ is a regular value of the moment map, we recover Witten's result \cite{Witten92} identifying the polynomial part of these integrals as the ordinary integral of the image of the class under the Kirwan map to the symplectic quotient.

\end{abstract}
 \maketitle    
     
\section*{Introduction}
Symplectic geometry (like many fields of geometry before it) has received a tremendous infusion of ideas from work on a single physics-inspired example.  The celebrated paper of Atiyah and Bott \cite{AB82} created a boom in symplectic geometry when they interpreted the moduli space of a Riemann surface as the symplectic reduction of the space of connections over that surface considered as a Hamiltonian space.  They were able to read off information about its cohomology by considering the square of the moment map as a Morse function (that this Morse function is exactly the Yang-Mills action reveals part of the physics inspiration behind their work).  Kirwan \cite{Kirwan84} built on these ideas to prove that for any compact Hamiltonian space the square  of the moment map is  an equivariantly perfect Morse function (or Kirwan-Morse function, a weaker but sufficient notion), and used this to give an algorithm for computing the Betti numbers of the symplectic quotient.   Not long afterwards inspiration came a second time from the same example: Witten \cite{Witten92} used quantum field theory ideas and some inventive symplectic geometry to find surprising formulas for the intersection pairings of the cohomology of the moduli space of a Riemann surface.  

Even for finite-dimensional Hamiltonian spaces Witten's techniques were not entirely rigorous, since he assumed regularities of the critical points of the squared moment map that do not typically hold. Jeffrey and Kirwan  were able to reproduce his key results both for moduli space \cite{JK98a} and for general compact Hamiltonian spaces \cite{JK95} by replacing his main technique (which he called nonabelian localization) with an older technique of Duistermaat-Heckman \cite{DH83} called abelian localization.  Specifically they were able to relate intersection pairings in the rational cohomology of a symplectic quotient to certain integrals of equivariant forms on the full Hamiltonian space.  In the intervening decade symplectic geometry has been a booming field with much of the work centering on the cohomology of the symplectic quotient and its relationship to the topology of the original Hamiltonian space (a sampling includes \cite{Kalkman95,GK96,Vergne96,BV97,LMTW98,MS99,Paradan00,Kiem04,TW03,JKKW03,BTW04}).

A remarkable feature of the work that has been driven by this one example is that none of the work actually applies to the original example (this is not entirely true:  as mentioned above Jeffrey and Kirwan \cite{JK98a} manage to prove Witten's formulas for moduli space, but the geometry of the space of connections which inspired these results is entirely circumvented).  The Hamiltonian space of interest is the space of connections, a space which is not just noncompact but in fact infinite dimensional,  and thus a far cry from the finite-dimensional compact Hamiltonian spaces on which we usually focus.  The failure of the results to apply to noncompact Hamiltonian spaces is particularly striking since the two most basic examples of Hamiltonian manifolds, $T^*G$ for $G$ a Lie group and any symplectic representation of a Lie group, are both noncompact.

Naturally one would like to extend Jeffrey and Kirwan's approach to the noncompact setting.  This may very well be possible, but one of the key benefits of reducing to the abelian case is the convexity results of Atiyah \cite{Atiyah82}, which apply only to compact spaces.  More precisely, Witten's results apply to Hamiltonian spaces for which $0$ is a regular value of the moment map.  In this situation a neighborhood of $\mu^{-1}(0)$ is always a Hamiltonian space with no fixed points for any subtorus, so any attempt to reduce questions of the topology of the reduced space to questions about the fixed points seems doomed.  In particular, a number of authors have extended Duistermaat-Heckman localization to noncompact settings \cite{PW94,Paradan00,Libine05}.  All these versions of Duistermaat-Heckman describe the induced measure on the Lie algebra at points away from $0,$ while in the case of a neighborhood of $\mu{-1}(0)$ when $0$ is regular  the integrals in question give measures on the Lie algebra with support entirely at $0.$ 

Witten's original approach to nonabelian localization, however, makes no apparent reliance on compactness.  In fact, if one were willing to join Witten in ignoring the analytic details and gave a sketchy introduction to equivariant cohomology in the Cartan model, Witten's argument could fit into a first year graduate course in differential geometry.   The intrinsic simplicity of his argument suggests that even in the compact case it may be illuminating and productive to work out the analytic details, assuming they are tractable.

They are indeed tractable, and this is the approach we take in this paper.  What makes them tractable, and in fact not very difficult, is the ability to avoid the central problem:  That the critical points of the square of the moment map, to which the integrals in question are supposed to localize, are in general singular spaces to which the differential geometry of forms and integration do not readily apply.  It seems likely that Witten's nonabelian localization has much more to tell us, but that to make further progress will require understanding these singular spaces.  for example, Paradan \cite{Paradan00} argues that the contribution to the Basic Integral from $\mu^{-1}(0)$ is still a polynomial even when $0$ is not a regular value.  This suggests that it still localizes to an integral of some sort of cohomology class over the (now singular) reduced space.  Since the higher critical sets can be built easily from $\mu^{-1}(0)$ of a related Hamiltonian space, we could then hope that the same is true for all critical sets.  

The paper is organized as follows.  Section 1 gives the local characterization of an arbitrary Hamiltonian space due to Guillemin and Sternberg \cite{GS84a,GS84b}.  It uses the local characterization to describe the critical points of the square of the moment map, and extend Kirwan's proof that the square of the moment map is an equivariantly perfect Morse-Kirwan function.  Section 2 reviews the Cartan model for equivariant cohomology, expressing the Cartan and Kirwan maps explicitly in this language. Section 3 defines the Basic Integral and computes key estimates for it.  Section 4 assumes that $0$ is a regular value of the moment map and proves the main theorem, that in this case the Basic Integral is the integral over the reduced space of the image under the Kirwan map of a certain form (which has polynomial dependence on $\epsilon$) plus additional contributions which are exponentially damped in $\epsilon.$  Thus cohomological integrals on the reduced space can be calculated by computing the Basic Integral over the full Hamiltonian space.

 \section{The  Local Structure Of Hamiltonian Spaces}\label{sc_local}

Let $M$ be a finite-dimensional smooth Hamiltonian space (not
necessarily compact):  That is a smooth manifold with symplectic form $\omega,$ acted on symplectomorphically by the
compact Lie group $G$ with Lie algebra $\lieg,$ and with moment map
$\mu \colon M \to \lieg^*.$ If $p \in M$ and $\xi \in \lieg$ we will
write $V\!\xi$ for the vector field associated to the infinitesimal
action of $\lieg$ on $M$ and $V_p\xi$ for the value of this vector
field at the point $p.$   Then the moment map condition is
\begin{equation} \label{eq:moment}
\omega \intprod V\!\phi = d\mu \dotprod \phi
\end{equation}
for all $\phi \in \lieg,$ where we use $\intprod$ to represent the interior product between a vector field and a form (or a tangent vector and a form at a point) with the convention that $v \intprod \omega = (-1)^{\deg(\omega)}\omega \intprod v.$   Choose an invariant inner product
$\bracket{\, \cdot \, , \, \cdot \,}_\lieg$ on
$\lieg.$  This inner product determines an identification $\star
\colon \lieg \to \lieg^*$ whose inverse we will also call $\star,$  so that
$\mu^\star\colon M \to \lieg.$  Finally, choose an almost
complex structure for $M$ compatible with the group action, that is to
say an invariant metric $\bracket{\, \cdot \, , \, \cdot \,}_M$ and an operator $J$ on
the tangent space such that $J^2=-1$ and $\bracket{x,y}=Jx\intprod
\omega\intprod y.$  
 
 \subsection{Local Characterization}
Guillemin and Sternberg (\cite{GS84a,GS84b})  give a local characterization of
a Hamiltonian space which will be crucial for what follows.    They
show that for any point $p \in M,$ the Hamiltonian  
space $M$ is determined in an equivariant neighborhood of $p$ by
the value of $\mu$ at $p,$ the Lie subgroup $H$ fixing $p$ and its Lie
subalgebra $\lieh,$  and the symplectic action of $H$ on the tangent space
$T_pM.$     More specifically, let $p \in M,$ with isotropy group $H \subset G,$ whose Lie algebra
   is $\lieh,$ and  define
   $\beta=\mu^\star_p$ and $K \subset G$ 
the stabilizer of $\beta,$ with $\liek$ its Lie algebra (so that
$\lieh \subset \liek$).  Let $Y$ be the subspace of $T_pM$ of vectors which are omega-orthogonal and orthogonal to $V\liek,$ the space of directions $V\!\phi$ for $\phi \in \liek.$  This is an $H$ representation and the symplectic form restricts to a symplectic form $\omega_X$ on $X.$  On the space
\[G \times ( \liek \oplus X)\]
 define the action of $G$ by the left action on the first component, and define the action of $H$ diagonally from the right action on the first component (applied to the inverse), the natural action on the second component, and the adjoint action on the third component.  Define a closed invariant two-form $\omega$ at the point $(g,\nu + x)$ by
\begin{equation}\label{eq:standard_omega}
\omega=\bracket{d\nu - \frac{1}{2} [\beta+\nu,g^{-1}dg],g^{-1}dg}_G + \frac{1}{2} dx\intprod \omega_X \intprod dx.
\end{equation}
and moment maps for the two actions 
\begin{equation}\label{eq:G_moment}
\mu_G  \phi = \bracket{\Ad_g(\beta + \nu),\phi}_G 
\end{equation}
\begin{equation}\label{eq:H_moment}
\mu_H  \phi = \bracket{\nu,\phi}_G  + \frac{1}{2}x \intprod \omega \intprod \phi x.
\end{equation}
The symplectic reduction by $\mu_H$ (i.e., the quotient of $\mu^{-1}(0)$ by the action of $H$ is  $G$-Hamiltonian space isomorphic to
\[G \times_H (\liek/\lieh \oplus X)\]
where we will interpret $\liek/\lieh$ as the subspace of $\liek$ perpendicular to $\lieh.$

One easily checks that the point $(1,0) \in G \times ( \liek \oplus X)$ is in $\mu_H^{-1}(0),$ has isotropy group $H,$ a tangent space isomorphic to $T_pM$ as an $H$-space, and moment value $\mu_G=\beta^\star.$  Therefore by
\cite{GS84b}[Thm. 41.2] there is an isomorphism of Hamiltonian spaces from $G \times_H (\liek/\lieh \oplus X)$ to a neighborhood of $p$ in $M.$  In the future we will refer to the choice of such an isomorphism as ``choosing a standard neighborhood of $p.$''

\subsection{The Square of the Moment Map}
  We are interested in the critical
sets of the nonnegative function $\bracket{\mu,\mu}_\lieg=|\mu|^2$ on $M.$

A point $ p\in M$ is a critical point for $|\mu|^2$ means that at $x,$
$d\bracket{\mu^\star,\mu^\star}_\lieg= 2\bracket{d\mu^\star,\mu^\star}_\lieg =d\mu
 \mu^\star= \omega \intprod V\!  \mu^\star =0.$  Since $\omega$ is
nondegenerate, to say that the one-form $V_p\mu^\star \intprod \omega$ is
zero at $p$ is to say that $V\!\mu^\star$ is zero at $p,$ and thus the
critical points of $|\mu|^2$ are exactly the zeros of the vector field
$V\!\mu^\star.$  Equivalently, critical points are the zeros of the one-form 
\begin{equation}\label{eq_lambda}
\lambda\intprod v \defequals \bracket{V\!\mu^\star, v}_M
\end{equation}
for $v$ a tangent vector on $M.$

If $p \in M$ and $G \times_H (\liek/\lieh \oplus X)$ is a standard neighborhood around $p$ then $p$ is a critical point for
the square of the moment map if and only if the $H$-orbit of $(1,0)$
 is critical in the standard neighborhood. This is equivalent to saying $V\!\beta=0,$ or
$\beta \in \lieh.$  
\begin{proposition} \label{pr_local_critical} Let $Z$ be the set of $x
\in X$
such that $\beta x=0$ and $Q_x=0,$ where $\bracket{Q_x,\phi}= \frac{1}{2} x \intprod
\omega_X \intprod \phi x$ for all  $\phi \in \lieh.$       The
connected component 
of the critical set  of  $|\mu|^2$ in $G \times_H(\liek/\lieh \oplus X)$  containing the $H$-orbit
of $(1,0)$  is the $H$-orbit 
of all points 
$(g, z)$ where $g \in G $ and $z \in Z$. 
This space is an algebraic variety.
\end{proposition}

\begin{pf}
Critical points of $|\mu|^2$ are points where $V\mu^\star=0.$  At a
point $(g, \nu + x)$ in $G \times (\liek/\lieh
\oplus X)$  we have
\[\mu \cdot \phi=\bracket{Ad_{a}(\beta + \nu),\phi} + \frac{1}{2} \phi x\intprod \omega_X \intprod x.\]
For a point $(g, \nu+x)$ in $\mu_H^{-1}(0)$ to descend to a point  for which $V_G\phi=0$ means  $V_G\phi$ is in $V_H \lieh.$  This
requires that $\Ad_{g^{-1}}\phi \in \lieh,$ 
$[\Ad_{g^{-1}}\phi,\nu]=0,$ and $\Ad_{g^{-1}}\phi\dotprod x=0.$
The first condition when $\phi=\mu^\star$ implies that $\nu=0.$  The second
implies nothing additional, and the third implies that
$(\beta+Q_x)x=0.$   Thus the critical points are in general those
for which $\nu=0$ and $(\beta+Q_x)x.$  The latter condition
implies that $\bracket{\beta+Q_x,Q_x}=0,$ which in turn implies that
$|\mu|^2 = |\beta+Q_x|^2 = |\beta|^2-|Q_x|^2.$  If there is a
path of  critical points connecting this to $(1,0),$ the value of
$|\mu|^2$ would be constant, which implies that $Q_x=0$ and hence
$\beta x=0.$ Of course
all points of this form are critical and are obviously path
connected to $(1,0),$ so the connected component of the critical set
includes these points. We have only to show all other solutions are
separated from this set. 

All other solutions have $\beta  x \neq 0.$  Since $\beta$
commutes with $\lieh \subset \liek,$  
write $X$ as a sum of orthogonal and $\omega$-orthogonal symplectic
submodules $X_\alpha,$ on each of which $|\beta  x| =\beta_\alpha
|x|$ for some positive $\beta_\alpha.$  If $x$ satisfies $(\beta +
Q_x) x=0,$ but not $\beta x=0$ then breaking $x$ into its
components there is a nonzero $x_\alpha$ such that $-\beta 
x_\alpha = Q_x  x_\alpha \neq 0.$  Since $|\beta  x_\alpha|=
\beta_\alpha |x_\alpha|$  at every solution $|Q_x| \geq
\min_\alpha \beta_\alpha$ holds.  Thus $|Q_x| > \frac{2}{3} \min_\alpha
\beta_\alpha$ and $|Q_x| < \frac{1}{3} \min_\alpha
\beta_\alpha$ separate $Z$ from all other solutions.
\end{pf}

\begin{corollary} \label{cr:describe_critical}
The set of critical points of $|\mu|^2$ on the Hamiltonian space $M$ is a discrete union of closed connected components, each of which is locally an algebraic variety and on each one of which the value of $\mu$ lies in a single coadjoint orbit.
\end{corollary}

\begin{corollary}
If $\mu$ is proper (that is the inverse image of compact sets is compact) then $|\mu|^2$ is a minimally degenerate equivariantly perfect 
Morse function in the sense of Kirwan \cite{Kirwan84}. 
\end{corollary}
\begin{pf}

In order for $|\mu|^2$ to be minimally degenerate we need to show that the critical set is a discrete union of compact sets on each of which $|\mu|^2$ is constant, and that for each of the sets there is a locally closed submanifold $\Sigma$ containing the critical set as a minimum and at each point in the critical set the tangent space to $\Sigma$ is a maximal subspace of the full tangent space on which the Hessian of $|\mu|^2$ is positive semidefinite.  The description of the critical sets is exactly the content of the previous corollary, together with the properness of $\mu.$  The existence of such a $\Sigma$ follows by the argument given by Kirwan unmodified, as does the equivariant perfection of this function.
\end{pf}

\section{Equivariant de Rham Cohomology}

An excellent reference on equivariant cohomology is \cite{GS99}, which gives a more complete and sophisticated treatment of everything in sections 2.1 and 2.2.

\subsection{Equivariant Forms}
Let $\mathcal{P}(\lieg)$ be the (graded) algebra of all complex-valued polynomial functions of $\lieg,$
 $\mathcal{S}(\lieg)$ be the algebra of complex-valued Schwartz functions on $\lieg$ (that is, any combination of derivatives of the function times any power of $|\phi|$ approaches $0$ as $\phi\to \infty,$ with the supremums of these products as seminorms),  $\mathcal{D}(\lieg)$ be the space of complex-valued tempered distributions, which is to say continuous linear functionals on $\mathcal{S}(\lieg),$ and $\mathcal{F}(\lieg)$ be the space of all continuous linear functionals on $\mathcal{P}(\lieg).$ Here and in the sequel we represent functions on
$\lieg$ as formulas in a dummy variable $\phi \in \lieg.$ Each of the function spaces ($\mathcal{P}(\lieg),$ $\mathcal{S}(\lieg)$) is an algebra and $\mathcal{P}(\lieg)$ acts by multiplication on $\mathcal{S}(\lieg),$ inducing various actions of the function spaces on the dual spaces ($\mathcal{D}(\lieg),$ $\mathcal{F}(\lieg)$) all represented by multiplication.  Also there are natural embeddings $\mathcal{P} \subset \mathcal{D},$ $\mathcal{S} \subset \mathcal{D},$ and $\mathcal{S} \subset \mathcal{F},$ sending $f(\phi)$ to $f(\phi) \dint \phi,$ where $\dint\phi$ represents Haar measure on $\lieg.$  By analogy with this embedding we will represent the pairing between a function space and its dual by $\int_\lieg \,\cdot.$

  If $M$ is a smooth
manifold and $\mathcal{X}$ represents one of $\mathcal{P}, \mathcal{S}, \mathcal{D}, \mathcal{F}$ we can define $\\Omega(M) \hat{\tensor}
\mathcal{X}(\lieg^*)$ to be smooth sections  of
the bundle over $M$ which at each point $p \in M$ 
is the tensor product of $\Lambda(T_pM) \tensor \mathcal{X}(\lieg).$  Here smooth means that when any element of the given space dual to $\mathcal{X}$ is paired with the second factor, the result is a smooth ordinary form.  When $\mathcal{X}$ is $\mathcal{P},$ this is an
algebra graded by the form degree plus twice the polynomial degree.  If $G$ acts smoothly on $M$ then $G$ acts
naturally and consistently on these bundles (with the diagonal action of $G$ acting
naturally on forms and by the dual of the adjoint action on
functions on $\lieg$), so we may
speak of the the $G$-invariant elements of each space.  These are
respectively the  \emph{$\mathcal{X}$-equivariant forms  on $M,$} though when $\mathcal{X}$ is $\mathcal{P}$ we drop the $\mathcal{P}$ and simply say \emph{equivariant forms on $M.$}
 The exterior
derivative $d$ is defined on   
all four bundles, so consider the \emph{equivariant derivative}
\begin{equation}\label{eq_equivariant_derivative}
D \alpha= d\alpha + i V\!\phi\intprod \alpha
\end{equation}  
where  $V\!\phi$ represents the linear map from $\lieg$ to vector
fields on $M.$  Note that $D$ is an equivariant map which
 increases degree by one and
satisfies $D^2=0$ on invariant elements.  Thus there are four equivariant cohomologies $H^{*,\mathcal{X}}_G(M),$ where in the case $\mathcal{X}=\mathcal{P}$ (the only case where the cohomology has an integer grading, the others have only a $\ZZ/2$ grading) we drop the $\mathcal{P}$ and write $H^{*}_G(M),$ the equivariant cohomology of $M.$ Equivariant differential forms give a
model for the cohomology with complex coefficients of the homotopy quotient $M_G,$ which is the
geometric significance of everything we do in this paper, but which is
mentioned for the last time here.

We say an $\mathcal{X}$-equivariant form has compact support if the closure of the set of points in $M$ where $\alpha$ is a nonzero function on $\lieg$ is compact.  Equivariant $D$
preserves both these concepts and we call the cohomology generated by
compactly-supported  $\mathcal{X}$-equivariant
forms $H^{*,\mathcal{X}}_{G,\compact}(M)$. 
The various products among $\mathcal{P},$ $\mathcal{S},$ $\mathcal{D},$ and $\mathcal{F}$ extend to products on the various equivariant forms by wedging the form component.  For example if $\alpha$ is an $\mathcal{S}$-equivariant form and $\beta$ is a $\mathcal{D}$-equivariant form then $\alpha\beta$ is an $\mathcal{F}$-equivariant form.  Note that equivariant $D$ satisfies the Leibniz rule on all such products.  Finally, an $\mathcal{F}$-equivariant form can be paired with $1$ (i.e. integrated) to get an ordinary form, which can be integrated over an invariant submanifold $N$ (assuming $M$ is oriented, and that either the original form was compactly supported or $N$ is compact) by taking on the component of appropriate degree.  Because $V\!\phi \intprod$ lowers form degree, 
\[\int_N \int_\lieg D\alpha \dint \phi= \int_N \int_\lieg d\alpha \dint \phi= \int_{\partial N} \int_\lieg \alpha \dint \phi\]
which is zero when $N$ has no boundary.  
This fact can be viewed as an
equivariant version of Stokes theorem and when applied to $N=M$ descends for example to a well-defined pairing on cohomology
\[H^{*,\mathcal{S}}_{G}(M)\times H^{*,\mathcal{D}}_{G,\compact}(M) \to \CC\]
and likewise with the compact subscript on the other factor.

\subsection{The Cartan and Kirwan Maps}
If the group action is locally free, the homotopy quotient retracts to the ordinary quotient and thus the equivariant cohomology is isomorphic to the ordinary cohomology of the quotient.  This isomorphism can be made completely explicit on the level of equivariant forms.

Let $P \to N$ be an orbifold principal $G$-bundle, which is to say locally $P$ can be identified with $G \times_H V$ where $H$ is a finite subgroup of $G$ and $V$ is an $H$-module, so that the $G$ orbit of each point in $G \times_H V$ is a fiber of the map $P \to N.$  Let $A$ be  a connection for this bundle,  i.e. an equivariant $\lieg$-valued one-form on $P$ such that $A\intprod V\!\phi=\phi$ for all $\phi \in \lieg.$ Let $P_A$ be the operator on $TP$ which sends a tangent vector $v$ to its projection onto the $A=0$ subspace, $P_Av=v - V\!A\intprod v.$  If $\alpha$ is a form on $P,$ define $P_A^*\alpha$ so that $v\intprod P_A^*\alpha = P_A^*(P_A(v) \intprod \alpha),$ i.e. $P_A^*\alpha$ is $\alpha$ projected onto the subspace of forms zero on all vertical vectors.  This map extends naturally to equivariant forms.   Define the \emph{Cartan map}
$\Cartan \colon \Omega(P) \hat{\tensor}
\mathcal{P}(\lieg) \to \Omega(P)$  by
\begin{equation}\label{weyl_map}
\Cartan(\alpha(\phi))=P_A^*(\alpha(iF_A))
\end{equation}
where $F_A$ refers to the $\lieg$-valued curvature two-form of the connection and its placement in parentheses denotes substituting its value  for $\phi$ in the second tensor factor of $\alpha,$ thus producing a form to be wedged with the first tensor factor.
\begin{proposition}
The Cartan map descends to a grade-preserving isomorphism from the complex of equivariant forms  to that of basic (i.e. horizontal and invariant) forms on $P$  inverting the natural imbedding.   Composing with the natural isomorphism of the complex of basic forms on $P$ with ordinary forms on $N,$ we get a map which descends to an isomorphism
\[\Cartan \colon H^*_G(P)\to H^*(N).\]
\end{proposition}

 \begin{pf}
 If $\alpha$ is an equivariant form on $P,$ it is clear that $\Cartan(\alpha)$ is invariant, by the equivariance of $P_A$ and $F_A.$  It is also clear that $\Cartan(\alpha)$ is horizontal, since $F_A$ is horizontal and the range of $P_A^*$ is horizontal.  Finally, it is clear that the Cartan map is an algebra homomorphism.  So for the homomorphism of complexes we need only show that the Cartan map intertwines the equivariant and ordinary derivatives, which can be checked locally.
 
 To do this consider a chart $V$ on which a finite subgroup $H$ of $G$ acts, and an equivariant isomorphism of $G \times_H V$ with a neighborhood in $P.$  $A$ defines an $H$-invariant one-form $\tau$ on $V$ with values in $\lieg,$ by $\tau_{v}\intprod \xi= A_{(1,v)}\intprod (0,\xi).$ A form on $G \times_H V$ is an $H$-invariant form on $G \times V.$  Since $G$ only acts on the first factor,  
 \[\left[\left(\Omega(G\times V)\right)^H \hat{\tensor} {\mathcal P}(\lieg^*)\right]^G \iso \left(\left[\Omega(G) \hat{\tensor} {\mathcal P}(\lieg^*) \right]^G\times \Omega(V)\right)^H\]  
 as complexes. 

Since $D$ and $d$ satisfy the Leibniz rule, we can check the intertwining on generators of the complex.  These are forms on $V,$ one-forms on $G$  $\bracket{\xi,g^{-1} dg }_G$ for $\xi \in \lieg,$ and functions $\bracket{\xi, g^{-1}\phi g}_G$ for  $\xi \in \lieg.$  That $D$ and $d$ are intertwined on the first set of generators is obvious.  For the second class
\begin{eqnarray*}
\Cartan(D(\bracket{\xi,g^{-1}dg })) &=& \Cartan(\bracket{\xi, g^{-1} dg g^{-1} dg } + i \bracket{\xi, g^{-1} \phi g})\\
&=& \bracket{\xi, \frac{1}{2}[\tau, \tau]} -\bracket{\xi, d\tau + \frac{1}{2}[\tau,\tau]}\\
&=& -\bracket{\xi,d\tau}\\
&=& -d(\bracket{\xi,\tau})\\
&=& d(\Cartan(\bracket{\xi,g^{-1}dg })).
\end{eqnarray*}
 For the third class
\begin{eqnarray*}
\Cartan(D(\bracket{\xi,g^{-1} \phi g }))&=& \Cartan(\bracket{\xi, [g^{-1}dg ,g^{-1} \phi g] } \\
&=& -i\bracket{\xi, [\tau, d\tau +\frac{1}{2} [\tau, \tau]]}\\
 &=&  i \bracket{\xi, [d\tau,\tau]}\\
&=& i d(\bracket{\xi,d\tau + \frac{1}{2}[\tau,\tau]})\\
&=& d(\Cartan(\bracket{\xi,g^{-1} \phi g })).
\end{eqnarray*}

Finally, to see that its inverse is the natural embedding of basic forms into equivariant forms, since it is the identity on basic forms, we need only check that every closed equivariant form is cohomologous to a basic form.  This requires defining certain operators on the complex of equivariant forms. 

We write $\frac{\partial}{\partial \phi}$ for the formal derivative with respect to $\phi,$ which we view as a function on $\lieg$ with values in equivariant forms.  Thus the operator 
\[\Phi=A\cdot \frac{\partial }{\partial \phi}\]
 denotes (viewing the connection $A$ as a form tensored with a Lie algebra element) applying this operator on the equivariant form to the second tensor factor and wedging the first tensor factor with the result.  By a similar logic $VA \intprod $ applies $V$ to the second tensor factor to get a tangent vector, takes the interior product with the form on which the operator acts to get a new form, and wedges the first factor with the result.  Now a straightforward calculation shows
\[ D \Phi + \Phi D = 
dA \cdot \frac{\partial}{\partial \phi}    + i(  \phi \cdot \frac{\partial}{\partial \phi} + V\!A \intprod \,)\]
where the two new operators in the above expression are defined similarly.  The first operator in the parentheses ($\phi \cdot \frac{\partial}{\partial \phi}$) multiplies any homogenous polynomial by its degree, and thus gives a grading of the space of equivariant forms into eigenvalues.  Similarly the second operator in the parentheses ($V\!A \intprod\,$) grades the space into eigenspaces with nonnegative integral eigenvalues, representing the ``number of form degrees in vertical directions.''  Since the two commute, they give a grading by their sum, call it the total degree, such that the total degree zero piece consists of basic forms on $P.$  Notice that the term not in parentheses ($dA \cdot \frac{\partial}{\partial \phi}$) strictly lowers total degree.   Thus if $\alpha$ is a closed equivariant form whose maximum total degree piece has degree $p>0,$ then  $\alpha + \frac{i}{p}D(\Phi\alpha)$ has strictly lower degree, and thus by induction $\alpha$ is cohomologous to a total degree zero form.
  \end{pf}

 Now suppose that $M$ is a Hamiltonian space with a proper moment map, and that $0$ is a regular value of $M,$ i.e. that $d\mu$ is onto for all points with $\mu=0.$

\begin{proposition}\label{pr:orbifold}
 If $d\mu$ is onto for each point of $\mu^{-1}(0),$ then $G$ acts on $\mu^{-1}(0)$ with finite stabilizers. In this case $\mu^{-1}(0)$ is a smooth manifold and an orbifold principal bundle over the quotient $M_{\text{red}}=\mu^{-1}(0)/G,$  which has an orbifold symplectic structure $\omega_0.$  \end{proposition}

\begin{pf}
If $d\mu$ is onto at some point, then by the moment map condition $V\!\phi$ is nonzero for all $\phi \in \lieg,$ so that the isotropy group must be finite.

If the isotropy group is finite at some $z$ with $\mu_z=0,$ then a standard neighborhood looks like
\[G \times_H (\lieg \oplus X)\]
where $H$ is a finite subgroup and $X$ is a symplectic vector space on which $H$ acts.  The subspace on which $\mu=0$ is $G \times_H X.$  The image of this space in the quotient by $G$ is isomorphic to $X/H,$ the quotient of a vector space by a finite-dimensional group action.     If another standard neighborhood $G \times_K X'$ contains $z,$ we argue the diffeomorphism of standard neighborhoods lifts to a local diffeomorphism of $X$ and $X'.$  This guarantees that a covering collection of standard neighborhoods form an orbifold atlas for $M_{\text{red}}=\mu^{-1}(0)/G.$

To see this, we can assume by equivariance that the standard neighborhood $G \times_K X'$ is chosen so that $z$ is the image of a point $(1,x).$  Then $H$ is the subgroup of $K$ which fixes $x,$ so that a neighborhood of $x$ is a representation $V$ of $H,$ and does not intersect with its image under any other elements of $K.$  Then $G \times_K X'$ is diffeomorphic locally to $G \times_H V,$ and this induces a diffeomorphism between $V$ and $X.$    

A symplectic structure on an orbifold is a choice of $H$-invariant symplectic form on $V$ for each chart $(V,H),$ which is preserved by the overlap maps.  Clearly $\omega_X$ is an invariant form on each vector space $X,$ and it is immediate that it is preserved by the overlap map.
\end{pf}

 The imbedding of $\mu^{-1}(0)$ into $M$ gives a map of equivariant cohomology which when composed with the Cartan map gives the \emph{Kirwan map} $\Kirwan \colon H^*_G(M) \to H^*(M_{\text{red}}).$  The fact that the square of the moment map is equivariantly perfect means that this map is surjective.

 \section{Equivariant Integration and Localization}
 
 For this section let $M$ be a Hamiltonian space with a proper moment map, and $\epsilon$ be a positive real parameter.

\subsection{Localization}
The moment map condition guarantees that the equivariant form 
\[\omega + i \mu \dotprod \phi\]
is closed, and thus represents an element of $H^*_G(M).$  Here the exponentiation is interpreted as its power series.  Suppose now that $\alpha$ is an  equivariant form on $M$, so that 
\[\alpha \exp(\omega + i \mu \dotprod \phi - \frac{\epsilon}{2} |\phi|^2)\]
is an $\mathcal{S}$-equivariant form which is closed and/or compactly-supported if $\alpha$ is.  

On the other hand consider an invariant ordinary one-form $\lambda$ on $M$ (which is therefore also an equivariant one-form).

\begin{lemma}\label{lm:Dlambda}
For each nonnegative real $t$
\[\int_0^t \exp(sD\lambda) \dint s \lambda\]
gives a $\mathcal{D}$-equivariant form satisfying $D\left(\int_0^t \exp(sD\lambda) \dint s \lambda\right)= \exp(tD\lambda)-1.$    Thus $\exp(tD\lambda)$
is a closed $\mathcal{D}$-equivariant form which is $\mathcal{D}$-cohomologous to $1.$
Further, on a submanifold of $M$ on which $\lambda\intprod V\!\phi$ is never the zero functional on $\phi,$ the limit of this integral as $t$ approaches infinity exists in the $\mathcal{D}$-topology and satisfies   $D\left(\int_0^\infty \exp(sD\lambda) \dint s \lambda\right)= -1.$
\end{lemma}

\begin{pf}
 We interpret the exponential and the integral in terms of power series,  and at a point in $M$ write $\lambda \intprod V\!\phi$ as $\bracket{\xi,\phi}$ for some $\xi\in \lieg,$ so that the integral is a sum of terms of the form
\[\bracket{\text{FORM}} \int_0^t s^k \exp(i s \bracket{\xi, \phi}) \dint s \]
which as a functional on some test function $f(\phi) \in \mathcal{S}(\lieg)$ is
\[\bracket{\text{FORM}}\int_0^t s^k \int_\lieg  f(\phi) \exp(is \bracket{\xi, \phi}) \dint \phi \dint s = \int_0^t s^k \widehat{f}(s\xi) \dint s \]
where $\widehat{f}$ is the Fourier transform of $f$ (ignoring arbitrary constants) and thus is well-defined.  So $\int_0^t \exp(sD\lambda) \dint s \lambda$ is a $\mathcal{D}$-equivariant form whose equivariant derivative is 
\[\int_0^t \exp(sD\lambda) D\lambda \dint s = 1-\exp(tD\lambda).\]

If $\lambda \intprod V\!\phi$ is never zero then $\xi\in \lieg$ as defined in the previous paragraph is never zero, so we get
\[\int_0^\infty s^k \widehat{f}(s\xi) \dint s\]
which converges since $\widehat{f}$ is Schwartz.

\end{pf}

\subsection{The Basic Integral}
Since $\mu$ is proper by Corollary \ref{cr:describe_critical}  identify the critical values of $|\mu|^2$ as
\[0 \leq  r_1 < r_2 < \cdots\]
(the sequence may be finite or infinite) and as long as $r\in \RR^+$ is regular, i.e. satisfies $r \neq r_i$ $\forall i \in \NN$ then 
\[M_r \defequals \{p \in M \,|\, |\mu_p|^2 \leq r\}\]
is a compact manifold with compact boundary.  
 Recall that the symplectic form gives a natural orientation to $M$ and hence $M_r$ and thus  integration over $M_r$ when $r$ is a regular value of $|\mu|^2$ is well-defined.   

Let $\lambda$ be the invariant one-form on $M$ which for any tangent vector $v$ gives
\begin{equation}\label{eq:lambda_def}
\lambda \intprod v = \bracket{V\!\mu^\star, v}.
\end{equation}
The $\lambda\intprod V\!\phi$ is zero exactly when $V\!\mu^\star$ is zero, which in turn happens exactly at the critical points of $|\mu|^2.$

For a equivariant form $\alpha,$ for any nonnegative real number $t$ and for any regular value $r$ of $|\mu|^2$ define the \emph{Basic Integral}
\begin{equation}\label{eq:regularized_int}
\BI(\alpha, r, t) \defequals \frac{1}{K}\int_\lieg   \int_{M_r} \alpha \exp[\omega + i \mu \dotprod \phi -\frac{\epsilon}{2} |\phi|^2 + t D\lambda] \dint \phi
\end{equation}
where $\lambda$ defined in Equation (\ref{eq:lambda_def}) and
\begin{equation}\label{eq:int_factor}
K= \vol(G) (2\pi)^{\dim(G)}.
\end{equation}

The following estimates are crucial to the calculations that follow.

\begin{lemma}\label{lm:bound_higher}
Suppose $\alpha$ is an equivariant form and $r$ and $s$  are regular values of $|\mu|^2$ with $s<r.$ Then
\[|\BI(\alpha,r,0)-\BI(\alpha,s,0)| < 
\text{POLYNOMIAL}(\epsilon^{\pm 1/2})\exp(-\frac{s}{2\epsilon})\]
where the coefficients of the polynomial depend on $r.$   In other words the contribution to the Basic Integral at $t=0$ of points with large values of $|\mu|$ is exponentially damped.
\end{lemma}

\begin{pf}
\begin{eqnarray*}
& &\left|\frac{1}{K}\int_\lieg  \int_{M_r-M_s} \alpha(\phi)  \exp(\omega + i \mu \dotprod \phi - \frac{\epsilon}{2} |\phi|^2)\dint \phi\right| \\
&=&  \frac{1}{K} \left|  \int_{M_r-M_s}  \exp(\omega)   \int_\lieg  \alpha(\phi)  \exp(i \mu \dotprod \phi - \frac{\epsilon}{2} |\phi|^2)   \dint \phi \right|\\
&=&  \frac{1}{K} \left|  \int_{M_r-M_s}  \exp(\omega)  \exp(-\frac{1}{2\epsilon} |\mu|^2) \int_\lieg  \alpha(\phi+i \mu^\star/\epsilon)   \exp(-\frac{\epsilon}{2} |\phi|^2)  \dint \phi 
 \right|\\
&=&  \frac{1}{K} \left|  \int_{M_r-M_s}  \exp(\omega)  \exp(-\frac{1}{2\epsilon} |\mu|^2)
\text{POLYNOMIAL}(\epsilon^{\pm1/2})
 \right|\\
&\leq& \exp(-\frac{s}{2\epsilon})\left|  \text{POLYNOMIAL}(\epsilon^{\pm1/2})\right|.
\end{eqnarray*}
Here the coefficients of the polynomial can be bounded by certain integrals over $M_r.$
\end{pf}

\begin{lemma}\label{lm:large_t}
Suppose that $\alpha$ is a \emph{closed} equivariant form and that $r$ is a regular value of $|\mu|^2.$ Then
\[ \lim_{t\to \infty} \BI(\alpha, r, t)\]
exists and differs from $\BI(\alpha, r, 0)$ by
\[ \text{POLYNOMIAL}(\epsilon^{\pm 1/2})\exp(-\frac{C}{2\epsilon})\]
where the coefficients of the polynomial and $C$ depend on $r.$
\end{lemma}

\begin{pf} Suppose $t_1<t_2 \in \RR.$
\begin{eqnarray*}
&&|\BI(\alpha,r, t_2)-\BI(\alpha, r, t_1)|\\
&=&\frac{1}{K}\left|  \int _\lieg  \int_{M_r} \alpha(\phi) 
\exp(\omega + i \mu \dotprod \phi -\frac{\epsilon}{2} |\phi|^2) 
(\exp(t_2D\lambda)-\exp( t_1 D\lambda) ) \dint \phi \right| \\
&=& \frac{1}{K}\left|  \int_\lieg   \int_{M_r} \alpha(\phi) 
\exp(\omega + i \mu \dotprod \phi -\frac{\epsilon}{2} |\phi|^2) 
D\left(\int_{t_1}^{t_2} \exp(sD\lambda) \lambda \dint s\right) \dint \phi\right| \\
&=&\frac{1}{K}\left|  \int_\lieg   \int_{\partial {M_r}} \alpha(\phi) 
\exp(\omega + i \mu \dotprod \phi -\frac{\epsilon}{2} |\phi|^2) 
\int_{t_1}^{t_2} \exp(sD\lambda) \lambda \dint s \dint \phi \right| \\
&=& \frac{1}{K}\Big|  \int_{\partial {M_r}} \int_{t_1}^{t_2}  \exp(\omega  +sd\lambda)\lambda
  \int_\lieg \alpha(\phi)
\exp(i \bracket{\mu^\star + sV^\star V\!\mu^\star, \phi} -\frac{\epsilon}{2} |\phi|^2) 
 \dint \phi \dint s\Big|. 
\end{eqnarray*}
 Completing the square yields
 \begin{eqnarray*}
&=& \frac{1}{K}\Big| \int_{\partial {M_r}}  \int_{t_1}^{t_2}     \exp(\omega  +sd\lambda)\lambda
 \int_\lieg  \alpha(\phi + \frac{1}{2\epsilon}(\mu^\star + sV^\star V\!\mu^\star) 
\exp(-\frac{\epsilon}{2} |\phi|^2)\dint \phi \\
&&\qquad \cdot \exp( -\frac{1}{2\epsilon} |\mu|^2  -\frac{s}{2\epsilon} |V\!\mu^\star|^2  -\frac{s^2}{2\epsilon} |V^\star V\!\mu^\star|^2) \dint s 
 \Big| \\
&=& \frac{1}{K}\Big| \int_{\partial {M_r}} \int_{t_1}^{t_2}     \exp(\omega  +sd\lambda)\lambda
\text{POLYNOMIAL}(\epsilon^{\pm1/2},s) \\
&&\qquad \cdot \exp( -\frac{1}{2\epsilon} |\mu|^2  -\frac{s}{2\epsilon} |V\!\mu^\star|^2  -\frac{s^2}{2\epsilon} |V^\star V\!\mu^\star|^2) \dint s
 \Big| \\
&\leq &   \frac{1}{K}\Big| \int_{\partial {M_r}}    \exp( -\frac{1}{2\epsilon} |\mu|^2)  
\text{FORM} \int_{t_1}^{t_2}    
\text{POLYNOMIAL}(\epsilon^{\pm1/2},s)  \\
&&\qquad  \cdot \exp(  -\frac{s^2}{2\epsilon} |V^\star V\!\mu^\star|^2) \dint s 
 \Big|.
\end{eqnarray*}
 For a fixed $\epsilon,$ since the $s$ integral is a polynomial times a Gaussian, this quantity is bounded by $\exp(-t_1^2 \min(|V^\star V\!\mu^\star|^2)/(2\epsilon)),$ where $  |V^\star V\!\mu^\star|^2$ is bounded below since $r$ is a regular value.  The limit of the difference can be written as a telescoping sum of such differences, which decrease hypergeometrically and hence the sum converges.    On the other hand choosing $t_1=0$ we see that there is a polynomial times $\epsilon^{\pm 1/2}$ which times $\exp(-C/(2\epsilon))$ bounds the difference regardless of $\epsilon$ or $t_2.$
\end{pf}

\subsection{The Basic Integral as a Sum of Contributions} The large $t$ limit of the Basic Integral is a sum of contributions from the critical points of $|\mu|^2,$ as is illustrated in the following.

\begin{lemma} \label{lm:r_independence}
If $r$ and $s$ are regular values of $|\mu|^2$ with no critical values between them and $\alpha$ is a closed equivariant form then
\[\lim_{t\to \infty} \BI(\alpha, r, t)=\lim_{t\to \infty} \BI(\alpha, s, t).\]
\end{lemma}

\begin{pf}
This follows directly from Lemma \ref{lm:Dlambda}.
\end{pf}

\begin{corollary} \label{cr:contribution}
For each $i$ choose $r_i'$ and $r_i''$ such that $r_{i-1}<r_i' < r_i < r_i'' < r_{i+1}.$  Define $r_1'=0$ and if $r_i$ is the maximum critical value choose any $r_i''>r_i.$  Then given a closed  equivariant form $\alpha$ the quantity
\[C_i(\alpha)=\lim_{t\to \infty} \BI(\alpha, r_i'',t) - \BI(\alpha, r_i'',t)\]
exists and is independent of the choice of $r_i'$ and $r_i''.$  Further, for any regular value $r$ of $|\mu|^2$ 
\[\lim_{t\to \infty} \BI(\alpha, r, t)= \sum_{r_i<r} C_i.\]
In other words the large $t$ limit of the Basic Integral up to $r$ is the sum of the contributions from each  critical set below $r.$  The contribution $C_i(\alpha)$ when $r_i>0$ is bounded by 
\[\text{POLYNOMIAL}(\epsilon^{\pm 1/2} \exp(-\frac{r_i-\delta}{2\epsilon})\]
where $\delta$ can be made as small as we like.
\end{corollary}

\begin{lemma} \label{lm:lambda_invariance}
Let $\alpha$ be any closed equivariant form, let $r_i$ be a critical value of $|\mu|^2,$ let $N$ be a compact manifold with boundary containing a neighborhood of the critical set corresponding to $r_i$ and no other critical points of $|\mu|^2,$ and let $\lambda'$ be the result of an isotopy of $\lambda$  such that the points of $M$ at which  $\lambda \colon V\!\phi$ is the zero functional on $\lieg$ remain fixed through the isotopy.  Then
\[
C_i(\alpha) = \lim_{t\to \infty} \int_\lieg \int_N \alpha  \exp(\omega + i \mu \phi + t D\lambda' -\frac{\epsilon}{2}|\phi|^2) \dint \phi.
\]
\end{lemma}

\begin{pf} By Lemma \ref{lm:Dlambda} the limit above with $\lambda$ replacing $\lambda'$ is equal to $C_i(\alpha).$ Define $\lambda''$ to agree with $\lambda'$ in a neighborhood of the critical set  but to agree with $\lambda $ near the boundary of $N.$  Then 
\begin{eqnarray*}
&&\frac{1}{K}\int_\lieg  \int_{N} \alpha \exp(\omega + i \mu \dotprod \phi -\frac{\epsilon}{2} |\phi|^2) \\
&& \qquad \cdot \left(\exp(tD\lambda) - \exp(tD\lambda'')\right)\dint \phi\\
&=&
\frac{1}{K}\int_\lieg  \int_{N} D\Big(
\alpha \exp(\omega + i \mu \dotprod \phi -\frac{\epsilon}{2} |\phi|^2) \\
&&\qquad \cdot \left( \int_0^t \exp(sD\lambda) \dint s \lambda - \int_0^t \exp(sD\lambda'') \dint s \lambda'' \right)\Big)\dint \phi\\
&=&
\frac{1}{K}\int_\lieg  \int_{\partial N} 
\alpha \exp(\omega + i \mu \dotprod \phi -\frac{\epsilon}{2} |\phi|^2)\\
&&\qquad \cdot  \left( \int_0^t \exp(sD\lambda) \dint s \lambda - \int_0^t \exp(sD\lambda'') \dint s \lambda'' \right)\dint \phi\\
&=&
0
\end{eqnarray*}
so that $\lambda$ and $\lambda''$ give the same contribution.  On the other hand by replacing $N$ by a smaller neighborhood (again by Lemma \ref{lm:Dlambda}), we can assure that $\lambda'$ and $\lambda''$ agree on $N$ and thus give the same contribution.
\end{pf}

\begin{proposition} \label{pr:bound_higher}
Suppose that $\alpha$ is a  closed equivariant form and $r $ is a regular value of $|\mu|^2.$  Then the large $t$ limit of the Basic Integral (\ref{eq:regularized_int}) is equal to its contribution $C_0(\alpha)$ of the critical set with $\mu=0$ (as defined in Corollary \ref{cr:contribution}) plus a contribution bounded by $\exp(-c/\epsilon)$  for some $c.$    
\end{proposition}

\begin{pf}
This follows immediately from Lemma \ref{lm:bound_higher}.  
\end{pf}

\section{When Zero is a Regular Value of the Moment Map}

  The proof of the following result in the case of trivial isotropy group appears in \cite{GS84b}, the full statement appears in \cite{Jeffrey99}.  While the statement and proof are widely known to experts, to the author's knowledge no proof appears in the literature, so for the sake of completeness it is included here.

\begin{proposition}\label{pr:normal_form}
 If $0$ is a regular value of $\mu$ (i.e. if $d\mu$ is onto for each point of $\mu^{-1}(0)$) recall by Proposition \ref{pr:orbifold} the map $\pi\colon \mu^{-1}(0) \to M_{\text{red}}=\mu^{-1}(0)/G$ is a principal orbifold bundle and $M_{\text{red}}$ has an orbifold symplectic structure $\omega_0.$  Given a connection $A$ on this bundle, there is an isomorphism of Hamiltonian spaces between a neighborhood of $\mu^{-1}(0)$ in $M$ and the Hamiltonian space $\mu^{-1}(0) \times \lieg,$ with symplectic form and moment map at $(p,\nu) \in \mu^{-1}(0) \times \lieg$ given by 
\begin{equation}
\widetilde{\omega}=\pi^* \omega_0 + d\bracket{\nu,A}
\end{equation} 
\begin{equation}
\widetilde{\mu}=\nu^\star.
\end{equation}
\end{proposition}

\begin{pf} One readily checks that  $\widetilde{\omega}$  defines a closed  form which is nondegenerate at $\mu^{-1}(0),$ and therefore in a neighborhood.  Also $\widetilde{\omega}$  is manifestly $G$-invariant (with the diagonal action of $G$ on $\mu^{-1}(0) \times \lieg$) and satisfies the moment map condition with $\widetilde{\mu}.$  By Guilleman and Sternberg's local characterization \cite{GS84b}[Thm.41.2], it suffices to give an equivariant symplectic isomorphism between the zeros of the moment map in each case, and then extend it to an equivariant identification of the normal bundles which  preserves $d\mu.$

The equivariant symplectic isomorphism is of course the natural imbedding of ${\widetilde{\mu}}^{-1}(0)= \mu^{-1}(0) \times \{0\}$ into $M.$ Its equivariance is by naturality and it preserves $\omega$ by inspection.  Because  $d\mu$ is onto at every point it gives a trivialization of the normal bundle, identifying it with $\mu^{-1}(0) \times \lieg.$  This identification clearly is equivariant and takes $d\widetilde{\mu}$ to  $d\mu.$
\end{pf}

\begin{theorem}  \label{th:zero_contribution}
Suppose $\alpha$ is a closed equiviariant form, and  $0$ is a regular value for the moment map.  Then the contribution $C_0(\alpha)$ to the Basic Integral (\ref{eq:regularized_int}) from $\mu^{-1}(0)$ is
\[\int_\lieg \dint \phi \int_{M_{\text{red}}} \Kirwan(\alpha) exp(\omega_0 + \frac{\epsilon}{2} c_2)\]
where $\Kirwan$ is the Kirwan map, $M_{\text{red}}$ is the orbifold quotient $\mu^{-1}(0)/G$ and $c_2$ is the second Chern class of the bundle $\mu^{-1}(0) \to M_{\text{red}}.$  In particular it has polynomial dependence on $\epsilon.$
\end{theorem}

\begin{pf}
The contribution to the basic integral of $\mu^{-1}(0)$ is
\[\lim_{t\to \infty}\frac{1}{K}  \int_N \int_\lieg \alpha(\phi)\exp(\omega + i \mu \cdot \phi - \frac{\epsilon}{2}|\phi|^2 + t D\lambda) \dint \phi\]
where $N$ is a neighborhood of $\mu^{-1}(0)$ containing no other critical points in its closure.  By Proposition \ref{pr:normal_form}  we can take $N$ isomorphic to a neighborhood of $\mu^{-1}(0)$ in $\mu^{-1}(0) \times \lieg.$  The integral is unchanged if we replace $\alpha$ by something cohomologous, so using an equivariant homotopy we can replace $\alpha$ with a form that agrees with $\iota^*(\alpha)\times 1$ in a neighborhood of $\mu^{-1}(0)$ in $\mu^{-1}(0) \times \lieg$ ($\iota$ being the inclusion of $\mu^{-1}(0)$). By making $N$ sufficiently small this form agrees with $\iota^*(\alpha)\times 1$ (which we will abbreviate $\iota^*(\alpha)$) everywhere.  Thus
\[=\lim_{t\to \infty}\frac{1}{K} \int_{N \subset \mu^{-1}(0) \times \lieg} \int_\lieg \iota^*(\alpha)(\phi) \exp(\pi^* \omega_0 + d\bracket{\nu,A} + i \bracket{\nu,\phi} - \frac{\epsilon}{2}|\phi|^2 + t D\lambda) \dint \phi.\]
By Lemma \ref{lm:lambda_invariance}, we may isotope $\lambda$ provided the zeros of $\lambda\intprod V\!\phi$ do not change.  Since $\bracket{V\!\mu^\star, \,\cdot \,}_M$ and $\bracket{\nu, A}_G$ are both positive on the vector $V\!\mu^\star,$ interpolating between them linearly does not change the zeros.  Thus replacing $\lambda$ with $\bracket{\nu,A}$ does not change the limit, giving
\begin{eqnarray*}
&=&
\lim_{t\to \infty}\frac{1}{K} \int_{N \subset \mu^{-1}(0) \times \lieg} \int_\lieg \iota^*(\alpha)(\phi) \\
&& \qquad \cdot \exp(\pi^* \omega_0 + d\bracket{\nu,A} + td\bracket{\nu,A} + i \bracket{\nu,\phi} + it \bracket{\nu,\phi} - \frac{\epsilon}{2}|\phi|^2) \dint \phi.
\end{eqnarray*}
Notice that $\nu$ occurs throughout with the factor $(1+t)$ (because $\iota^*(\alpha)$ does not depend $\nu$) and thus we can rescale to eliminate $t$ except for the dependence of the region of integration.  In the large $t$ limit this becomes the integral over all $\lieg$
\begin{eqnarray*}
&=&
\frac{1}{K} \int_{\mu^{-1}(0) \times \lieg} \int_\lieg \iota^*(\alpha(\phi)) \\
&&\qquad \cdot \exp(\pi^* \omega_0 + d\bracket{\nu,A}  + i \bracket{\nu,\phi}  - \frac{\epsilon}{2}|\phi|^2) \dint \phi.
\end{eqnarray*}
Completing the square
\begin{eqnarray*}
&=&\frac{1}{K} \int_{\mu^{-1}(0) \times \lieg}  \left(\int_\lieg \iota^*(\alpha)(\phi+\frac{i}{\epsilon} \nu) 
 \exp(-\frac{\epsilon}{2}|\phi|^2) \dint \phi\right)\\
&&\qquad \cdot \exp(\pi^* \omega_0 + d\bracket{\nu,A}   - \frac{1}{2\epsilon}|\nu|^2)\\
&=&
\frac{1}{K} \int_{\mu^{-1}(0) \times \lieg}  \left(\int_\lieg \iota^*(\alpha)(\phi+  \frac{i}{\epsilon}  \nu) \exp(-\frac{\epsilon}{2}|\phi|^2) \dint \phi\right)\\
&&\qquad \cdot \exp(\pi^* \omega_0 + \bracket{d\nu, A} + \bracket{\nu,dA}   - \frac{1}{2\epsilon}|\nu|^2).
\end{eqnarray*}
Note the only occurrence of $d\nu$ is in $\exp(\bracket{d\nu,A}).$  Consider a basis of tangent vector at some point in $\mu^{-1}(0) \times \lieg$ which consists of an orthonormal basis of $\lieg,$ the image of this orthonormal basis under $V,$ and a basis of $A$-horizontal vectors in $\mu^{-1}(0).$  The top dimensional piece of this multiform is a sum of terms with $\bracket{d\nu,A}$ raised to various powers, but the only terms which are nonzero when applied to this basis are those where $\bracket{d\nu,A}$ is raised to $\dim(G),$ and on those terms the value is unchanged if $P_A^*$ is applied to all other forms in the product.  Thus
\begin{eqnarray*}
&=&\frac{1}{K} \int_{\mu^{-1}(0) \times \lieg}  \left(\int_\lieg P_A^*\circ \iota^*(\alpha)(\phi+ \frac{i}{\epsilon}  \nu) 
\exp(-\frac{\epsilon}{2}|\phi|^2) \dint \phi\right)\\
&&\qquad \cdot \exp(\pi^* \omega_0  + \bracket{\nu,P_A^*(dA)}   - \frac{1}{2\epsilon}|\nu|^2)\bracket{d\nu,A}^{\dim(G)}\\
&=&
\frac{1}{K} \int_{\mu^{-1}(0) \times \lieg}   \left(\int_\lieg P_A^*\circ \iota^*(\alpha)(\phi+ \frac{i}{\epsilon}  \nu)
 \exp(-\frac{\epsilon}{2}|\phi|^2) \dint \phi\right)\\
&&\qquad \cdot \exp(\pi^* \omega_0  + \bracket{\nu,F_A}   - \frac{1}{2\epsilon}|\nu|^2)\bracket{d\nu,A}^{\dim(G)}\\
&=&
\frac{1}{K} \int_{\mu^{-1}(0) \times \lieg}   \left(\int_\lieg P_A^*\circ \iota^*(\alpha)(\phi+ \frac{i}{\epsilon} \nu +i F_A) 
\exp(-\frac{\epsilon}{2}|\phi|^2- \frac{1}{2\epsilon}|\nu|^2) \dint \phi\right)\\
&&\qquad \cdot \exp(\pi^* \omega_0  + \frac{\epsilon}{2} |F_A|^2 ) \bracket{d\nu,A}^{\dim(G)}\\
&=&\frac{\vol(G)}{K}    \int_{\mu^{-1}(0)/G} \left(\int_\lieg \int_\lieg P_A^*\circ \iota^*(\alpha)(\phi+ \frac{i}{\epsilon} \nu +i F_A)
 \exp(-\frac{\epsilon}{2}|\phi|^2 - \frac{1}{2\epsilon}|\nu|^2) \dint \phi \dint \nu\right)\\
&&\qquad \cdot \exp( \omega_0  + \frac{\epsilon}{2} |F_A|^2 ) 
\end{eqnarray*}
where we have completed the square on $\nu$ and integrated the result over the vertical fibers, noting that the integral is constant in these directions and that the measure $ \bracket{d\nu,A}^{\dim(G)}$ is equal to Haar measure on the vertical fiber times Lebesgue  measure $\dint \nu$ on $\nu .$  Now changing the $\nu$ and $\phi$ variables to a single complex variable $z=\sqrt{\epsilon} \phi + i \nu/\sqrt{\epsilon}$ and noting that the integral of any complex polynomial against a complex Gaussian measure gives its constant term yields
\begin{eqnarray*}
&=&\frac{\vol(G)}{K}    \int_{\mu^{-1}(0)/G} \left(\int_{\lieg + i \lieg}  P_A^*\circ \iota^*(\alpha)(\frac{1}{\sqrt{\epsilon}}z +i F_A) \exp(-|z|^2) \dint z\right)\exp( \omega_0  + \frac{\epsilon}{2} |F_A|^2 ) \\
&=&\frac{\vol(G)(2 \pi)^{\dim(G)})}{K}    \int_{\mu^{-1}(0)/G}   P_A^*\circ \iota^*(\alpha)(i F_A) \exp( \omega_0  + \frac{\epsilon}{2} |F_A|^2 ) \\
&=&\frac{\vol(G)(2 \pi)^{\dim(G)})}{K}    \int_{\mu^{-1}(0)/G}   \Kirwan(\alpha) \exp( \omega_0  + \frac{\epsilon}{2} |F_A|^2 ) 
\end{eqnarray*}
\end{pf}

\begin{corollary}\label{cr:reduction_of_integral}
If $\alpha$ is a closed equivariant form, $r$ is a regular value of $|\mu|^2,$ $\mu$ is proper and $0$ is a regular value of $\mu$   then the Basic Integral $\BI(\alpha, r, 0)$ can be written uniquely as a sum of a polynomial in $\epsilon$ plus a term bounded by $\exp(-c/\epsilon)$ for some $c>0,$ the polynomial piece representing the contribution from $\mu^{-1}(0)$ as in Theorem \ref{th:zero_contribution}.  
\end{corollary}

\begin{pf}
We know that $\BI(\alpha, r,  0)$ differs from the large $t$ limit $\lim_{t\to \infty} \BI(\alpha, r, t)$ by a quantity bounded by $\exp(-c/\epsilon)$ for some $c>0$ by Lemma \ref{lm:large_t}.   On the other hand the large $t$ limit is a sum of contributions from $r=0$ and higher critical values by Corollary \ref{cr:contribution}.  The former is a polynomial in $\epsilon$ by Theorem \ref{th:zero_contribution}, the latter is bounded by $\exp(-c/\epsilon)$ for some $c>0$ by Proposition \ref{pr:bound_higher}.  Since a function can only be written in one way as a polynomial plus a term bounded by $\exp(-c/\epsilon),$ the result follows.
\end{pf}

In general we have no reason to believe that the integral over all of $M,$ that is the large $r$ limit of $\BI(\alpha, r, t)$ exists for a fixed $t$ or the large $t$ limit.  however, if it exists and converges sufficiently rapidly, the same results as above apply.  For example
\begin{proposition}
Suppose that $M$ is a Hamiltonian space with proper moment map and $0$ is a regular value of $\mu.$  Suppose also the symplectic volume of $M_r$ as a function of $r$ is such that $|\partial \Vol(M_r)/\partial r | < \exp(c \sqrt{r})$ for some $c>0.$   Suppose also that for some almost complex structure the supremum over all of $M_r$ of the norm of $\alpha$ (the norm as an ordinary form at each point times the norm as a symmetric tensor in $\lieg^*$) is also bounded by $\exp(c \sqrt{r}).$  Then 
\[\lim_{r\to \infty} \int_{M_r} \int_\lieg \alpha \exp(\omega + i \mu \phi -\frac{\epsilon}{2}|\phi|^2)\dint \phi\]
exists and is of the form a polynomial in $\epsilon$ plus a term exponentially damped in $\epsilon,$ the polynomial 
\end{proposition}

\begin{pf} Fix a regular value $r_0$ of $|\mu|^2.$
\begin{eqnarray*}
&&\left| \lim_{r\to \infty}
\BI(\alpha, r, 0)- \BI(\alpha, r_0, 0)\right|\\
&=& \frac{1}{K} \int_{M-M_{r_0}} \int_{\lieg} \alpha(\phi) \exp(\omega + i \mu \phi -\frac{\epsilon}{2}|\phi|^2)\dint \phi\\
&=& \frac{1}{K} \left|\int_{M-M_{r_0}} \exp(\omega  - \frac{1}{2\epsilon}|\mu|^2)\int_{\lieg} \alpha(\phi +_ i \mu^\star/\epsilon) \exp( -\frac{\epsilon}{2}|\phi|^2)\dint \phi\right|\\
&\leq& \frac{1}{K} \left|\int_{r_0}^\infty \exp(-\frac{r}{2\epsilon} + c \sqrt{r}) \text{POLY}(\sqrt{r}, \epsilon^{\pm 1/2}) \partial \Vol(M_r)/\partial r \dint r\right|\\
&\leq& \frac{1}{K} \left|\int_{r_0}^\infty \exp(-\frac{r}{2\epsilon} + 2c \sqrt{r}) \text{POLY}(\sqrt{r}, \epsilon^{\pm 1/2}) \dint r\right|\\
&\leq&  \exp(-\frac{k}{\epsilon}) \text{POLY}( \epsilon^{\pm 1/2}) 
\end{eqnarray*}
for some positive constant $k.$  Applying this to Corollary \ref{cr:reduction_of_integral} gives the result.

\end{pf}

\bibliographystyle{alpha}
\def\cprime{$'$}

\end{document}